# Proof of Descartes circle formula and its generalization clarified

*Jerzy Kocik*

Department of Mathematics, Southern Illinois University, Carbondale, IL 62901
jkocik@siu.edu

In his talk "Integral Apollonian disk Packings" Peter Sarnak asked if there is a "proof from the Book" of the Descartes theorem on circles [Sar]. A candidate for such a proof is presented in this note.

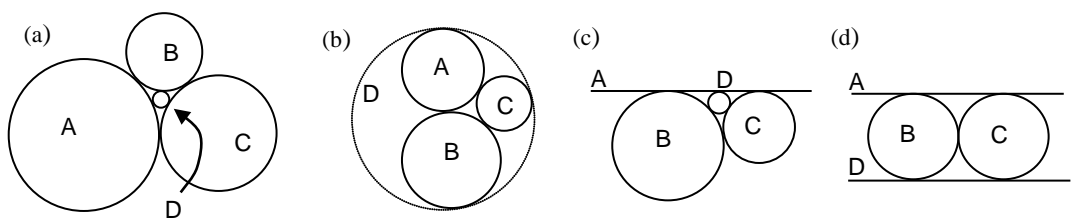

**Figure 1**: Four disks in a Descartes configuration — special cases Then the outer circle in (b) represents the boundary of an unbounded disk *outside* of circle D, for which we assume a negative radius and curvature. Each pair of disks are tangent *externally*.

**Theorem 1 (Descartes formula, 1643 [Ped]).** The curvatures $a$, $b$, $c$ and $d$ of four pair-wise externally tangent disks satisfy

$$(a + b + c + d)^2 = 2(a^2 + b^2 + c^2 + d^2). \qquad (1)$$

A system of such four pair-wise tangent disks is said to form **Descartes' configuration**.

What follows is a proof of a generalized theorem for disks in general configuration, which will be stated at the end, and from which (1) follows as a special case. This exposition is a succinct and clarified version of the proof in [JK].

**Proof:** Define a "scalar product" for the disks in the plane as the cosine of the angle under which the circles intersect. In case they do not intersect, we may extend the definition by the Law of Cosines:

$$\langle C_1, C_2 \rangle = \begin{cases} (d^2 - r_1^2 - r_2^2)/2r_i r_j & \text{in general} \\ \cos \varphi, & \text{if } C_1 \text{ and } C_2 \text{ intersect} \end{cases} \qquad (2)$$

where $d$ is distance between the centers, and $r_1$, $r_2$ are the corresponding radii of the disks, and $\varphi$ is the angle under which the disks intersect.

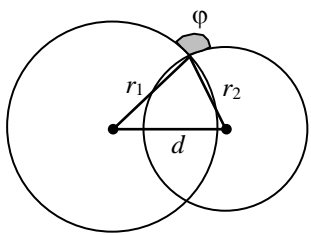

**Figure 2:** Intersecting disks/circles

Let $\mathbf{M} \cong \mathbb{R}^{3,1}$ denote Minkowski space with the metric (quadratic form) defined

$$g = \begin{bmatrix} -1 & 0 & 0 & 0 \\ 0 & -1 & 0 & 0 \\ 0 & 0 & 0 & 1/2 \\ 0 & 0 & 1/2 & 0 \end{bmatrix} \qquad (3)$$

There is a map that sends disks (circles) in the plane to space-like unit vectors in $\mathbf{M}$, namely, a disk centered at $(x, y)$ and radius $r$ is mapped to vector

$$\text{Disk}((x,y); r) \longrightarrow C = \begin{bmatrix} \dot{x} \\ \dot{y} \\ \beta \\ \gamma \end{bmatrix} = \begin{bmatrix} x/r \\ y/r \\ 1/r \\ (x^2 + y^2 - r^2)/r \end{bmatrix} \qquad (4)$$

where $\beta = 1/r$ is the curvature, $\dot{x} = x/r$, and $\dot{y} = y/r$ are "reduced coordinates.", Now, the inner product (2) becomes

$$\langle C_1, C_2 \rangle = C_1^\mathrm{T} g\, C_2 \qquad (5)$$

Indeed, with simple algebraic manipulations we get:

$$\begin{aligned}
\langle C_1, C_2 \rangle &= \frac{d^2 - r_1^2 - r_2^2}{2 r_1 r_2} = \frac{(x_2 - x_1)^2 + (y_2 - y_1)^2 - r_1^2 - r_2^2}{2 r_1 r_2} \\
&= -\frac{x_1 x_2}{r_1 r_2} - \frac{y_1 y_2}{r_1 r_2} + \frac{1}{2 r_1} \frac{x_2^2 + y_2^2 - r_2^2}{r_2} + \frac{1}{2 r_2} \frac{x_1^2 + y_1^2 - r_1^2}{r_1} \\
&= -\dot{x}_1 \dot{x}_2 - \dot{y}_1 \dot{y}_2 + \tfrac{1}{2}\beta_1 \gamma_2 + \tfrac{1}{2}\beta_2 \gamma_1
\end{aligned}$$

which has manifestly the form of an inner product defined by (3). In particular, the circle vectors are normalized: $\langle C, C \rangle = -\dot{x}^2 - \dot{y}^2 + \beta \gamma = -1$. Now, build a matrix by juxtaposing the four vectors representing the disks:

$$D = [C_1 | C_2 | C_3 | C_4] \qquad (6)$$

The collective scalar product of disks in (6) results in the Gramian $f$ defined by

$$f = D^\mathrm{T} g\, D \qquad (7)$$

It may be called "**configuration matrix**". For instance, for Descartes' configuration, i.e., four mutually tangent disks, the product of each pair of distinct disks is equal to 1 (since $\cos 0 = 1$):

$$\langle C_i, C_j \rangle = \begin{cases} 1 & \text{if } i = j \\ -1 & \text{if } i \neq j \end{cases} \qquad (8)$$

Take the inverse of both sides of (7) to get

$$F = (D^{-1})\, G\, (D^\mathrm{T})^{-1}$$

where $G = g^{-1}$ is the inverse of $g$, and $F = f^{-1}$ is the inverse of $f$. Next, multiply both sides of the equation by $D$ on the left and $D^\mathrm{T}$ on the right to arrive at

$$DFD^\mathrm{T} = G. \qquad (9)$$

And this is Descartes' theorem generalized! It may be formulated in the following "stand alone" form:

**Theorem (Generalized disk/circle formula):** Four disks in a general configuration in the plane satisfy the following matrix equation

$$D^T F D = G \qquad (10)$$

where $F = f^{-1}$ is the inverse of the "configuration matrix" $f$ the entries of which are defined by

$$f_{ij} = \langle C_i, C_j \rangle = \frac{d^2 - r_1^2 - r_2^2}{2 r_1 r_2}.$$

In particular, for the vector of four curvatures, $B = [\,a,\,b,\,c,\,d\,]^T$, formula $B^T F B = 0$ replaces Descartes formula (1).

**Example:** In Descartes configuration, Equation (8) turns (9) into:

$$\underbrace{\begin{bmatrix} -1 & 1 & 1 & 1 \\ 1 & -1 & 1 & 1 \\ 1 & 1 & -1 & 1 \\ 1 & 1 & 1 & -1 \end{bmatrix}}_{f} = \underbrace{\begin{bmatrix} \dot{x}_1 & \dot{y}_1 & \beta_1 & \gamma_1 \\ \dot{x}_2 & \dot{y}_2 & \beta_2 & \gamma_2 \\ \dot{x}_3 & \dot{y}_3 & \beta_3 & \gamma_3 \\ \dot{x}_4 & \dot{y}_4 & \beta_4 & \gamma_4 \end{bmatrix}}_{D^T} \underbrace{\begin{bmatrix} -1 & 0 & 0 & 0 \\ 0 & -1 & 0 & 0 \\ 0 & 0 & 0 & 1/2 \\ 0 & 0 & 1/2 & 0 \end{bmatrix}}_{g} \underbrace{\begin{bmatrix} \dot{x}_1 & \dot{x}_2 & \dot{x}_3 & \dot{x}_4 \\ \dot{y}_1 & \dot{y}_2 & \dot{y}_3 & \dot{y}_4 \\ \beta_1 & \beta_2 & \beta_3 & \beta_4 \\ \gamma_1 & \gamma_2 & \gamma_3 & \gamma_4 \end{bmatrix}}_{D} \qquad (11)$$

It is easy to see that $f^2 = 4I$, hence $F = f^{-1} = (1/4)f$. Scaling both sides by 4 we arrive at the extended Descartes' formula for curvatures, positions and co-curvatures (cf. [LMW]):

$$\underbrace{\begin{bmatrix} \dot{x}_1 & \dot{x}_2 & \dot{x}_3 & \dot{x}_4 \\ \dot{y}_1 & \dot{y}_2 & \dot{y}_3 & \dot{y}_4 \\ \beta_1 & \beta_2 & \beta_3 & \beta_4 \\ \gamma_1 & \gamma_2 & \gamma_3 & \gamma_4 \end{bmatrix}}_{D} \underbrace{\begin{bmatrix} -1 & 1 & 1 & 1 \\ 1 & -1 & 1 & 1 \\ 1 & 1 & -1 & 1 \\ 1 & 1 & 1 & -1 \end{bmatrix}}_{\tfrac{1}{4}F} \underbrace{\begin{bmatrix} \dot{x}_1 & \dot{y}_1 & \boxed{\beta_1} & \gamma_1 \\ \dot{x}_2 & \dot{y}_2 & \beta_2 & \gamma_2 \\ \dot{x}_3 & \dot{y}_3 & \beta_3 & \gamma_3 \\ \dot{x}_4 & \dot{y}_4 & \beta_4 & \gamma_4 \end{bmatrix}}_{D^T} = \underbrace{\begin{bmatrix} -4 & 0 & 0 & 0 \\ 0 & -4 & 0 & 0 \\ 0 & 0 & \boxed{0} & 8 \\ 0 & 0 & 8 & 0 \end{bmatrix}}_{4G} \qquad (12)$$

The part enclosed by rectangles is Descartes' formula (1).

**Remark 1:** The literature often presents Descartes theorem as a statement about circles. It is less confusing to view it as a theorem on configuration of non-overlapping **disks**. In particular, the "outer circle" in Figure 1(b) represents the boundary of an unbounded disk *outside* of circle $D$, for which we assume a negative radius and curvature. The straight lines in Fig. 1(d) represent halfplanes – disks of curvature 0, which are also tangent externally with the other disks.

**Remark 2:** The proof generalizes to $(n-1)$-spheres in $n$-dimensional space with Minkowski space $\mathbf{R}^{n+1,1}$ in place of $\mathbf{R}^{3,1}$ [JK].